\documentclass{article}
\usepackage{graphicx}
\usepackage{amsmath}
\usepackage{amssymb}
\usepackage{simplemargins}
\usepackage[hidelinks]{hyperref}

\settopmargin{2cm}
\setbottommargin{2cm}
\setleftmargin{2cm}
\setrightmargin{2cm}

\title{\large \textbf{Endless love: On the termination of a playground number game}}
\author{\small Iain G. Johnston \\ \vspace{-0.2cm} \footnotesize School of Biosciences, University of Birmingham, UK \\ \vspace{-0.2cm} \footnotesize i.johnston.1@bham.ac.uk}
\date{}

\begin{document}

\maketitle
\noindent A simple and popular childhood game, `LOVES' or the `Love Calculator', involves an iterated rule applied to a string of digits and gives rise to surprisingly rich behaviour. Traditionally, players' names are used to set the initial conditions for an instance of the game: its behaviour for an exhaustive set of pairings of popular UK childrens' names, and for more general initial conditions, is examined. Convergence to a fixed outcome (the desired result) is not guaranteed, even for some plausible first name pairings. No pairs of top-50 common first names exhibit non-convergence, suggesting that it is rare in the playground; however, including surnames makes non-convergence more likely due to higher letter counts (for example, `Reese Witherspoon LOVES Calvin Harris'). Different game keywords (including from different languages) are also considered. An estimate for non-convergence propensity is derived: if the sum $m$ of digits in a string of length $w$ obeys $m > 18/(3/2)^{w-4}$, convergence is less likely. Pairs of top UK names with pairs of `O's and several `L's (for example, Chloe and Joseph, or Brooke and Scarlett) often attain high scores. When considering individual names playing with a range of partners, those with no `LOVES' letters score lowest, and names with intermediate (not simply the highest) letter counts often perform best, with Connor and Evie averaging the highest scores when played with other UK top names.

\section*{Introduction}
The `LOVES' game, variously known as the `Love Calculator', `crush compatibility test' and others, is played in various guises in playgrounds and classrooms across the world \cite{roud2010lore, ref1}. The game has given rise to websites and apps peppering the web (for example, \url{http://www.lovecalculator.com}); and its whimsical nature has been cited as motivation in developing more rigorous game-theoretical approaches to assess compatibility \cite{bever2013love}. `LOVES' is a word-and-number game which attempts to assign a two-digit number to two character strings, by applying an iterated rule to a initial set of numbers, derived from the occurrence counts of certain letters within the strings. This note aims to explore the outcomes of this iterated rule for different real-world and general initial conditions, and for different game structures.

%

\subsection*{The game}
\begin{figure}
\centering
\includegraphics[width=0.6\linewidth]{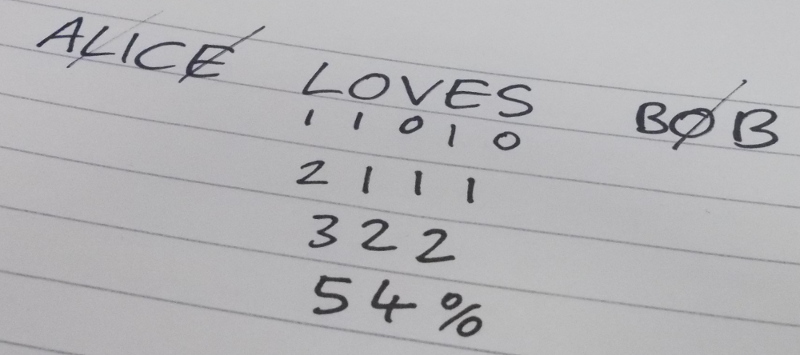}
\caption{The `LOVES' game played by Alice and Bob, as described in the text.}
\label{fig0}
\end{figure}

We begin with two names. Occurrence counts of each letter in the word `LOVES' in the two names are recorded, and are expressed as a string of integers, which we will call $s^{(0)}$. For example, if Alice were playing with Bob, we would count 1 `L', 1 `O', 0 `V's, 1 `E' and 0 `S's, so $s^{(0)} = \{1, 1, 0, 1, 0\}$ (Fig. \ref{fig0}).

We then move through the digits in the string from left to right. In this order, we sum each digit with the neighbour on its right, and append the sum to the rightmost position of new string $s^{(1)}$ (which starts off empty). So for Alice and Bob, the first pair of digits in $s^{(0)}$ ($1$ and $1$) sum to $2$, the next pair ($1$ and $0$) sum to $1$, the next to $1$, and the final pair also to $1$, so we arrive at the new string $s^{(1)} = \{2, 1, 1, 1\}$ (Fig. \ref{fig0}).

When we reach the end of our previous string, we repeat this process for the new string that has been formed, leading to another new string. This process is iterated until we reach a string containing fewer than three digits. Following this process, we see $s^{(2)} = \{2, 1, 1, 1\} \rightarrow s^{(3)} = \{3, 2, 2\} \rightarrow s^{(4)} = \{5, 4\}$, whereupon the iteration terminates, having produced a string of length $< 3$ (Fig. \ref{fig0}). The traditional conclusion of the game is to represent the two subjects' love as a percentage made from the final string: `Alice loves Bob 54\%!'. 

More generally, given an initial string of digits $s^{(0)}$, the game consists of applying an iterative rule to the current string $s^{(t)}$ to produce a new string $s^{({t+1})}$:

\begin{equation}
s^{(t+1)} = (s^{(t)}_1 + s^{(t)}_2)\,||\,(s^{(t)}_2 + s^{(t)}_3)\,||\,...\,||\,(s^{(t)}_{n_t-1} + s^{(t)}_{n_t}), \label{eqn1}
\end{equation}

where $s_i$ denotes the $i$th digit in string $s$, $n_t$ is the number of digits in string $s^{(t)}$, and $||$ is the concatenation operator.

Several variations exist, including changes to which letters or sets of letters are counted. Sometimes an additional modular number game, `FLAMES', is played to decide the mode of interaction: (Friendship, Love, and others).

\section*{Results I: Types of love}

\subsection*{General behaviour}
We observe three types of behaviour produced by different initial strings: \emph{termination}, \emph{looping}, and \emph{divergence}. Termination is illustrated by the example above and in Fig. \ref{fig0}: a string of length $< 3$ is reached after some finite number of iterations. The other two `endless love' cases never reach this target. Looping involves a periodic cycling between two or more different strings of length $\geq 3$, repeated forever. Divergence involves an unbound increase in string length as iterations continue. Table \ref{tab1} contains examples of these different behaviours.

\begin{table*}
\footnotesize
\begin{tabular}{| l | l |}
\hline
Termination & $\{1, 1, 0, 1, 0\} \rightarrow \{ 2, 1, 1, 1\} \rightarrow \{3, 2, 2\} \rightarrow \{5, 4\}$ \\
Looping & $\{1, 4, 1, 2, 1\} \rightarrow \{5, 5, 3, 3\} \rightarrow \{1, 0, 8, 6\} \rightarrow \{1, 8, 1, 4\} \rightarrow \{9, 9, 5\} \rightarrow \{1, 8, 1, 4\} \rightarrow ...$ \\
Divergence & $\{6, 6, 6, 6, 6\} \rightarrow \{1, 2, 1, 2, 1, 2, 1, 2\} \rightarrow \{3, 3, 3, 3, 3, 3, 3\} \rightarrow \{6, 6, 6, 6, 6, 6\} \rightarrow \{1, 2, 1, 2, 1, 2, 1, 2, 1, 2\} \rightarrow ...$ \\
\hline
\end{tabular}
\caption{Examples of the three types of outcome of the game.} 
\label{tab1}
\end{table*}

Two looping `motifs' can be identified, namely a period-2 cycle $\{9, 9, ...\} \rightarrow \{1, 8, 1, ...\}$ and a period-3 cycle $\{3, 3, 3, 3, ...\} \rightarrow \{6, 6, 6, ...\} \rightarrow \{1, 2, 1, 2, 1, ...\}$, where the ellipses in the first step given refer to at least one other digit. Depending on the subsequent digits, these motifs may also form part of divergent strings (and usually do if there is more than one other digit; hence, these motifs will usually not give rise to looping behaviour if included as part of an initial length-5 string).



\subsection*{Distribution of outcomes by starting set}
We first investigate the distribution of behaviours across different sets of initial strings. We will restrict our analysis\footnote{The extension to initial counts over 10 will be seen to fall naturally into our analysis when we consider longer starting strings.} to starting strings where all letter counts are under 10. Initially, we consider sets of strings classified by an integer label $c$. A set with label $c$ contains all strings of five integers which all have values less than $c$. Hence, $c=10$ corresponds to the set $\{0,0,0,0,0\} ... \{9,9,9,9,9\}$ of all five-digit strings; $c=2$ corresponds to the set $\{0,0,0,0,0\} ... \{1,1,1,1,1\}$ of all binary strings of length five. 

\begin{figure}
\centering
\includegraphics[width=0.6\linewidth]{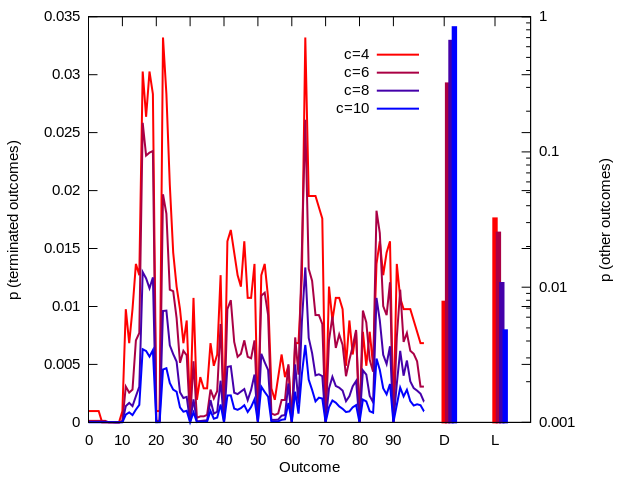}
\caption{Frequencies of different outcomes of `LOVES' when played with all different starting strings consisting of five digits less than $c$. Left hand side (linear vertical axis) shows individual results ($0\% ... 99\%$) of terminating strings; right hand side (logarithmic vertical axis) shows `endless love' cases: D -- divergent; L -- looping.}
\label{fig1}
\end{figure}

Fig. \ref{fig1} shows the distribution of game outcomes over sets with $c = 4, 6, 8, 10$. We can see that the proportion of divergent strings increases towards 1 as $c$ increases, suggesting that initial strings containing higher digits are more likely to diverge. Looping behaviour is rather more rare: the proportion of looping strings is under $0.03$ for $c=4$ and decreases as $c$ increases. The proportion of terminating strings decreases with $c$ in concert with the increasing likelihood of divergence. Of those strings that do terminate, the pattern of final outcomes displays complicated structure, which is to some extent conserved across different $c$ values. Multiples of 10 are notably rarer than other outcomes, as only one string can give rise to these final values: the outcome $10n$ is only produced by the string $\{n, 0, 0, 0, 0\}$. Regions of high occurrence occur across $c$ values in the ranges $15-22$, $50-55$, $62-68$ and $85-92$, providing a rather mixed outcome for potential lovers: the reason for this structure currently remains elusive. 

\subsection*{Distribution of outcomes by starting magnitude}
Noting that divergence appears to be more common for strings containing higher digits, we assign a descriptor of a string, its \emph{magnitude} $m$. For a string $s$ of length $L$ we define $m$ simply as the sum of all digits in the string $m = \sum_{i=1}^L s_i$. The observations above suggest that divergence propensity will increase with $m$.

\begin{figure}
\centering
\includegraphics[width=0.6\linewidth]{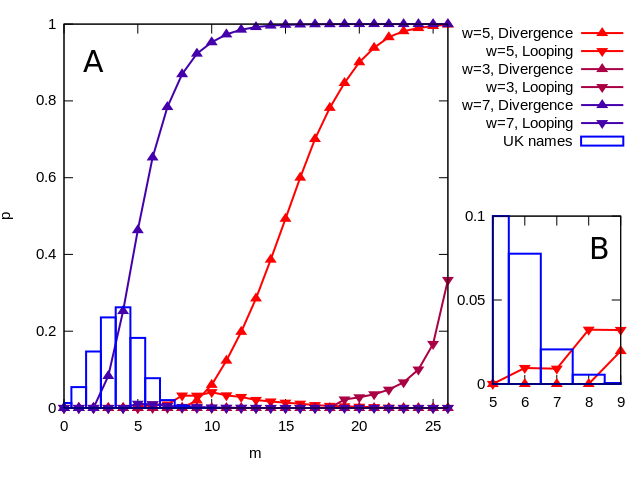}
\caption{Outcomes of games over the set $c=10$ as a function of start string magnitude $m$ and length $w$ (`LOVES' corresponds to $w=5$; other games are discussed later in the text). (A) Probability that a string diverges or loops with $m$; histogram shows the distribution of magnitudes found when playing `LOVES' with random pairings of popular UK first names (see text). (B) Same plot with focus on the region where UK name distribution overlaps with a nonzero looping probability for $w=5$.}
\label{fig2}
\end{figure}

To explore the influence of magnitude on the probability of `endless love', we subdivide the set $c=10$ into strings labelled by $m$, and compared the outcome behaviour of each of these $m$-labelled subsets. Fig. \ref{fig2} shows these results: there is a monotonic increase of divergence probability with $m$, with divergence probability saturating at 1 as $m$ increases. No strings with $m>27$ terminated (the terminating $m = 27$ case is illustrated by $\{7, 3, 9, 1, 7\}$, which terminates at $\{6, 5\}$ after 24 iterations).

Terminating strings with high $m$ tend to display motifs where several pairs of adjacent digits sum to 10. This property is clearly exhibited by the $m=27$ example above; others include $\{5, 2, 8, 2, 9\}$ ($m = 26$, $\{8, 5\}$ after 17 iterations) and $\{1, 5, 5, 5, 9\}$ ($m = 25$, $\{2, 2\}$ after 27 iterations). Intuitively, these strings immediately expand into strings containing $\{..., 1, 0, 1, 0, ...\}$ which then collapse, diminishing string length and yielding small numbers in the next iteration: a hint to players who want to ensure their love remains bound. In the Appendix we analyse in more detail the length of the `decay chains' by which an initial string evolves, and the behaviour of the string length as this process progresses.
 
\subsection*{Termination of combinations of real names}
The non-terminating behaviour of the game may be surprising to readers who did not encounter this behaviour in the playground \footnote{Although the looping behaviour has been previously been observed by some players; for example, by commenters on this discussion \url{http://www.lawoftheplayground.com/browse.php?type=subject&id=2449}}. We next seek to explain this by considering the game behaviour with initial strings that are likely to appear in real games. To this end, we use lists of the 50 most popular baby boys' and girls' first names in the UK in 2010 \cite{ref2} (thus presumably populating playgrounds today)\footnote{Readers may be interested that `Olivia' and `Oliver' occupy the female and male top spots; a selection of other names feature in Table \ref{highscoretable}.}. The initial number strings arising from all boy-girl, boy-boy, and girl-girl pairings of these names were computed, and the magnitudes of these `realistic' strings were compared to our findings in the previous section. There is indeed little overlap between the magnitudes arising from common UK name pairings and the magnitudes leading to non-termination (Fig. \ref{fig2} inset). Although a small amount of overlap with looping magnitudes did exist, no specific pairs of top-50 names exhibit this behaviour. It is straightforward to construct reasonable cases that do, however: consider `Steve-O LOVES Esmie' $\rightarrow \{0, 1, 1, 4, 2\}$ ($m = 7$), which collapses to the $\{9,9,1\} \leftrightarrow \{1, 8, 1, 0\}$ loop. Including surnames in the game makes achieving the magnitudes required for non-convergence easier: for example `Reese Witherspoon LOVES Calvin Harris' $\rightarrow \{1, 2, 1, 4, 3\}$ ($m = 11$) $\rightarrow \{9, 9, 6\} \leftrightarrow \{1, 8, 1, 5\}$.

\begin{figure}
\centering
\includegraphics[width=0.6\linewidth]{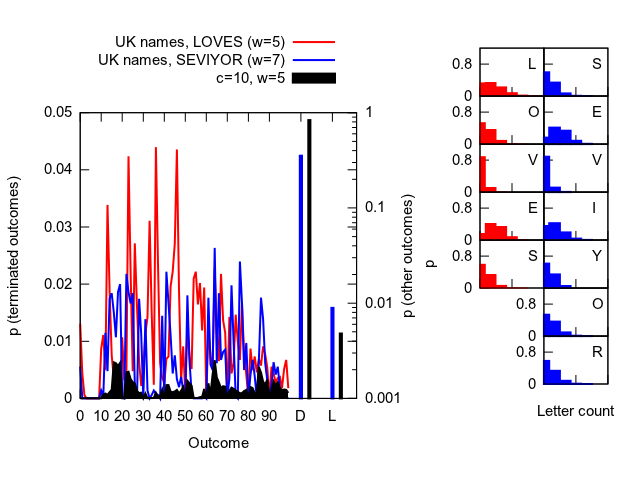}
\caption{(left) Outcomes of the `LOVES' and `SEVIYOR' games with the set of common UK name pairings, compared to `LOVES' with the set $c=10$ (as in Fig. \ref{fig1}). D -- divergent, L -- looping. (right) Distributions of letter counts among the set of common UK pairings.}
\label{fig2-5}
\end{figure}

The distribution of outcomes from the set of common UK names is shown in Fig. \ref{fig2-5}. There is substantially different structure compared to the full $c=10$ set, even when the higher probability of termination is accounted for: the landscape of peaks of high occurrence is more rugged and does not obviously correlate with the common regions previously noted, with a (rather unfortunate) skew of terminating strings towards lower final results. One reason for this is the dramatic difference in letter frequencies among UK names compared to a uniform sampling: V in particular is highly underrepresented (Fig. \ref{fig2-5} inset). Generally lower numbers in the central position of the start string may be expected to give rise to differences in game outcome. Our exhaustive search of $c=10$ led to the observation that for a count of zero `V's, compensatory high counts for other characters (two counts of 4, or one of 3 and one of 5), are required to prevent termination; the requirement for a compensatory high count even for a single V is illustrated in our example above. 

If we remove the infrequent V and replace it with commoner letters, we add a degree of nonsense to the literal interpretation of the game, but can discover looping behaviour among common UK names. For example, `Connor LORES Harrison' $\rightarrow \{0, 3, 3, 0, 1\}$, which collapses to the $\{9, 9, 4\} \leftrightarrow \{1, 8, 1, 3 \}$ cycle.

\subsection*{Other games}
We have seen that `endless love' is a rare occurrence in UK playgrounds. How much is this an unfortunate result of linguistics? A natural extension of the `LOVES' game is to consider different central words: varying the length of the initial string of numbers will intuitively impact on the outcome of the game. Here we consider the case where the central word consists of $w$ strictly different letters. Cases where the central word contains repeated letters (for example, `ELSKER' in Norwegian), and some elements in the initial number string are thus forced to be identical, constitutes a further extension to be explored in future.

As expected, we observe an increase in divergence propensity with increasing $w$, for the $c=10$ set (where the set of initial strings is now understood to consist of all strings of length $w$ consisting of numbers under $c$), illustrated in Fig. \ref{fig2}. The sigmoidal relationship between magnitude and divergence frequency is shifted to lower $m$ and sharpened as $w$ increases, indicating that lower magnitudes can achieve non-convergence.

An interesting result of increasing $w$ is the shift of divergent magnitudes into the range occupied by common UK name pairings. For example, if UK children were to play using the $w=7$ Turkish word `SEVIYOR', Fig. \ref{fig2} predicts that divergent behaviour would be much more common. This prediction is borne out by simulations pairing common UK names and playing `SEVIYOR' -- the proportion of divergent name pairings rises from zero to around a third of all pairings (Fig. \ref{fig2-5}), showing that `endless love' may be rather more common in non-English-speaking playgrounds.

\section*{Results II: The path to love}
\subsection*{Mean trajectories of strings}
To understand how players may predict their propensity for `endless' or bounded love, we now take a more analytic view of a given game. If we consider a chain as a `walk' in $(w, m)$ space, different pairs of adjacent numbers in the current string can give rise to two different walking behaviours. For a single pair ($w=2$), if the sum of the numbers $<10$, the effect of that pair on application of Eqn. \ref{eqn1} is to leave $m$ unchanged and reduce $w$ by 1. If the sum of the numbers $\geq 10$, $m$ is reduced by 9 and $w$ remains the same. 

A string can be viewed as a collection of $w-1$ pairs, where each digit except the first and last contribute to two pairs. The extension of the above $w=2$ picture is as follows. The maximum attainable magnitude in the next step is $m' = s_1 + s_w + 2\sum_{i=2}^{w-1} s_i$, and the maximum attainable length is $w' = 2(w-1)$. Each pair that sums to $< 10$ contributes $\Delta m = 0, \Delta w = -1$; each pair that sums to $\geq 10$ contributes $\Delta m = -9, \Delta w = 0$. It will readily be seen that the statistics of the string after an iteration of Eqn. \ref{eqn1} are $w \rightarrow w' + \Delta w, m \rightarrow m' + \Delta m$. 

Hence, for $\{9, 9, 9, 9, 9\}$, we have $m' = 72, w' = 8$; all four pairs sum to $> 10$, so we have $\Delta m = -36, \Delta w = 0$, giving $m = 36, w = 8$, the statistics of $\{1, 8, 1, 8, 1, 8, 1, 8\}$. As another example, for $\{1, 0, 9, 2, 2\}$ we have $m' = 25, w' = 8$; three pairs sum to $<10$ and one to $>10$ so we have $\Delta m = -9, \Delta w = -3$, giving $m = 16, w = 5$, the statistics of $\{1, 9, 1, 1, 4\}$. 

We write $n$ for the number of pairs in a string that sum to $\geq 10$. The total number of pairs is $w-1$: hence, $w-1-n$ pairs sum to $< 10$. The next iteration is then characterised by

\begin{eqnarray}
m & \rightarrow & m' - \Delta m = s_1 + s_w + 2\sum_{i=2}^{w-1} s_i - 9 n \label{mupdate} \\
w & \rightarrow & w' - \Delta w = 2(w-1) - (w - 1 - n) =  w - 1 + n. \label{wupdate}
\end{eqnarray}

We can begin to estimated the expected action of iterating Eqn. \ref{eqn1} by considering the dynamics that it provokes in $(w,m)$ space. To proceed, we record the average $(\Delta w, \Delta m)$ resulting from an application of Eqn. \ref{eqn1} to all strings that occupy a given point in $(w,m)$. This averaged quantity $(\langle \Delta w \rangle, \langle \Delta m \rangle)$ gives the expected movement in $(w,m)$ space for each point.

\begin{figure}
\centering
\includegraphics[width=0.8\linewidth]{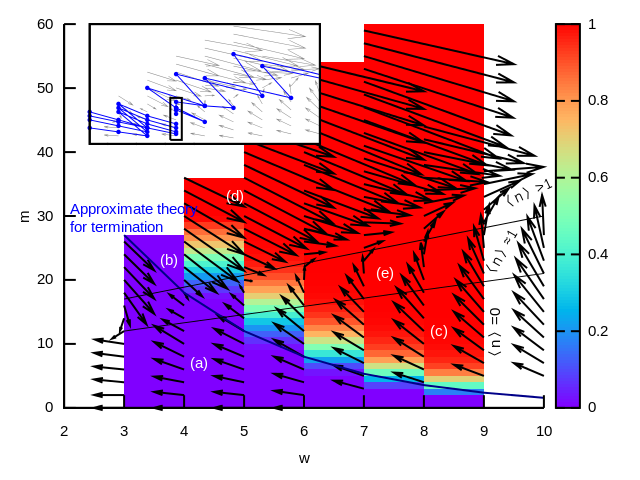}
\caption{Mean trajectories in $(w,m)$ space. Coloured background shows the probability of divergence for each $(w,m)$ point. Black vectors show step magnitude and direction, upon one iteration of the game rule, averaged over all strings corresponding to a given $(w,m)$. Black lines roughly divide $(w,m)$ space by $\langle n \rangle$, the number of length-increasing transitions. The blue line gives a theoretical estimate for the boundary of the region of divergence from Eqn. \ref{theoryline}. Inset shows some example trajectories from specific terminating (leftwards) and divergent (rightwards) strings, from different source strings within the black rectangle.}
\label{figsteps}
\end{figure}

Fig. \ref{figsteps} plots the averaged step $(\langle \Delta w \rangle, \langle \Delta m \rangle)$, taken over all strings with a given $(w,m)$. The structure of the expected step behaviour in the $(w,m)$ plane immediately sheds light on the dynamics we have seen so far. Low $(w,m)$ strings (Fig. \ref{figsteps}a) experience moderate increases in $m$ in concert with decreases in $w$. For sufficiently low starting $m$, the small magnitude of $m$ increase allows these dynamics to reach $w=2$ without exceeding a `critical line' of $m$. There is also a region of $w=3$ strings with high $m$ (Fig. \ref{figsteps}b) where an iteration produces an increase in $w$, but coupled with a sufficient loss of $m$ to keep strings in the low-$(w,m)$ region. 

Outside this region, two behaviours are visible. For high $w$ and low $m$ (Fig. \ref{figsteps}c), decreases in $w$ are induced by iteration, but accompanied by large increases in $m$. The trend to decreasing $w$ suggests that $n=0$ is common in this region (from Eqn. \ref{wupdate}). For high $m$ (Fig. \ref{figsteps}d), iterations induce decreases in $m$ accompanied by often dramatic increases in $w$, suggesting (Eqn. \ref{wupdate}) that $n>0$ in this region. The action of these two trends is to force strings towards a noticeable line of divergence around region (\ref{figsteps}e) and above, where iterations leads to increases in $w$ with moderate increases in $m$, with no upper bound. The slow changes in $w$ in this line of divergence suggests that $n \simeq 1$ in this region. This splitting of the phase portrait by $n$ is indeed what is observed, as can be seen in the demarcated regions of $\langle n \rangle$ ($n$ averaged over all strings for a given $(w,m)$) in Fig. \ref{figsteps}.

Fig. \ref{figsteps} of course only represents an `averaged' behaviour, but the agreement with the observed divergence propensities is striking. It can be seen that movement towards the $\langle n \rangle \simeq 1$ band in Fig. \ref{figsteps} is a central determinant of eventual behaviour. As long as $n = 0$, $w$ decreases; divergence requires steps with $n>0$. For a given $w$, the value of $m$ marking the transition from highly-likely termination to highly-likely divergence seems to correspond to the region where the averaged result of a step avoids the $\langle n \rangle \simeq 1$ band until $w \leq 4$. 

For $n = 0$, Eqn. \ref{mupdate} gives bounds on the step in magnitudes that occur with the $w \rightarrow w - 1$ transition. In the Appendix, we show that a prediction for the transition line separating convergent from non-convergent regions of $(w,m)$ space occurs approximately at

\begin{equation}
m^* \simeq \frac{18}{(3/2)^{w-4}}, \label{theoryline}
\end{equation}

which agrees with observations in Fig. \ref{figsteps}.


\subsection*{Termination result and magnitude}
\begin{figure}
\centering
\includegraphics[width=0.6\linewidth]{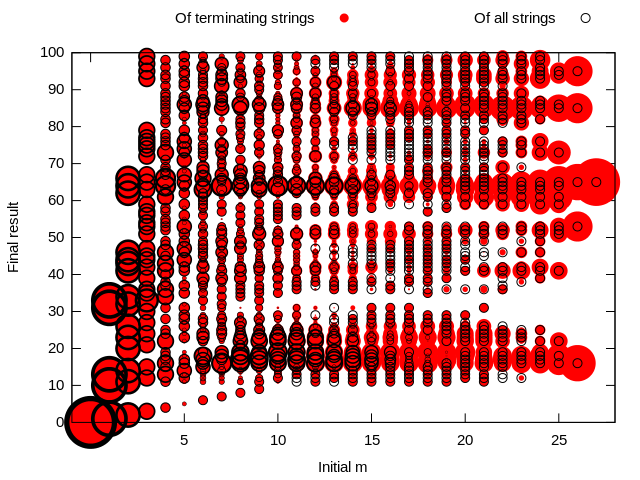}
\caption{Probability of a given final result $10 s_1 + s_2$ from a terminating string with initial magnitude $m$. Radii of filled circles are proportional to the log probability that a \emph{terminating} string of magnitude $m$ produces an output; radii of empty circles are proportional to the log probability that \emph{any} string of magnitude $m$ produces an output.}
\label{figoutcomes}
\end{figure}

The anticipated question, what properties of an initial string predict termination at a high final `love' outcome, remains challenging to answer. No simple properties of $w=5$ strings, including $m$, individual $s_i$ values, and linear combinations of $s_i$ values, and differences between adjacent $s_i$, were found to correlate with either termination magnitude or the final result at termination (defined as $10 s_1 + s_2$). Fig. \ref{figsteps} hints at a possible strategy for maximising final $m$: choosing an initial point in $(w,m)$ space that maximises $m$ while avoiding the $n=1$ zone. However, the nonlinearity in the iterative process, and the pronounced string-to-string variability in specific behaviours from any $(w,m)$, makes even this approach unreliable; and further extension from final magnitude to final result is also difficult. Fig. \ref{figoutcomes} shows the probability of different outcomes for different strings of magnitude $m$; while there is a slight trend for \emph{terminating} strings of higher $m$ to yield higher results than terminating strings of lower $m$, the decreasing number of terminating strings of high $m$ diminishes this trend when all strings are considered.

Even explicitly `backtracking' from a desired target string fails to yield much intuitive progress. In the Appendix we show that in the absence of length-expanding steps ($n=0$ everywhere), the initial string giving a final result $\{s_1, s_2\}$ is

\begin{eqnarray}
\left\{ a_3, a_2-a_3, a_1 - 2 a_2 + a_3, s_1 - 3 a_1 + 3 a_2 - a_3, s_2 - 3 s_1 + 6 a_1 - 4 a_2 + a_3 \right\},
\end{eqnarray}

for a set of digits $a_1, a_2, a_3$ that can to some extent be specified. However, the numerous inequalities that $a_1, a_2, a_3$ must satisfy to yield a valid decay chain (see Appendix) prevent the straightforward choice of a starting string to give a desired result; integer programming, which is NP-hard, is required to make progress within this picture \cite{integerprogramming}. The question of what intuitive features, if any, predict whether a given starting string will terminate with a high final value remains open.

One step towards an answer is to examine the set of initial strings that give rise to high final results. Many strings yield $99\%$, from the $m = 3$ $\{0,0,3,0,0\}$ to the $m=24$ $\{9,1,9,3,2\}$; of the $10^5$ strings with $c = 10$, 97 give $99\%$ and 731 give $95\%$ or above. However, many of these would not plausibly arise from real names; restricting strings to $m < 10$ and the number of `V's to be under 3 gives 76 strings giving $95\%$ or above. Only 10 actually occur in combinations from the list of top-50 UK names (Table \ref{highscoretable}). Heuristically, names with pairs of `O's or high counts of `L's are among the top scores, though this is certainly not a rigorous rule.

\begin{figure}
\centering
\includegraphics[width=0.6\linewidth]{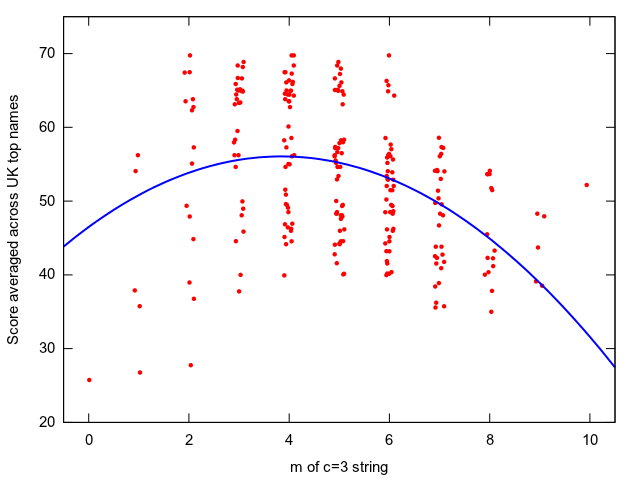}
\caption{Magnitude of individual $w=5$, $c=3$ strings against their averaged score when `LOVES' is played with all UK top names. For clarity, points are randomly jittered on the x-axis about their (integer) values. Line is a quadratic fit with $r = 0.41$.}
\label{figzeroes}
\end{figure}

One can also consider the individual names that score best when the game is played with other names. Of course, the performance of individual names is identical if they contain the same counts of the `LOVES' letters. We therefore consider the sets in $c=3$ (deemed realistic for real-world names) that have the highest average score when played with all 100 top UK names. Ranked by average score and omitting the string notation for brevity, the winners were 21012, 20110, 01102 and 00200, averaging $69.76\%$ across the set of names. The winning sets that themselves correspond to top-50 names were 02000 (Connor), scoring $67.42\%$ on average, and 00120 (Evie) scoring $66.69\%$. The bottom was the set 00000, scoring $25.75\%$, with corresponding names ranging from Adam to Ryan (and, sadly, the author's own name). One notable property of the set of high-scoring $c= 3$ sets is intermediate magnitudes: sets with $3 \leq m \leq 6$ score more highly on average when played with UK names than those with higher or lower $m$ (Fig. \ref{figzeroes}).

Ranked by number of $\geq 95\%$ scores, the winning set is 00222 (with 28 partners giving $\geq 95\%$), followed by 12010 and 21101, each with 26 partners. By this criterion, the highest-ranking sets with top-50 names were 02010 (Brooke) with 20 partners and 11100 (Olivia) with 16 partners; again, the 00000 set ranked lowest, tying with several others with 0 partners yielding $\geq 95\%$ scores. 

Players with names matching the high-scoring strings given here and in Table \ref{highscoretable} are more likely to score highly when playing `LOVES' with a randomly chosen partner. We note that `Vivian', while not a top-50 UK name, seems set to experience romantic success, matching the top-scoring 00200 string. 

\begin{table}
\scriptsize
\begin{tabular}{|| p{9cm} l l | l l ||}
\hline
\hline
\textbf{High scoring ($\geq 95\%$) complete starting strings} & & & & \\
\hline
\hline
Names scoring 99\%& & & Sets & Score \\
\hline
Joseph and Leo, Chloe, Charlotte & & & 02100, 22010 & 95 \\
Brooke and Lewis, Samuel, Scarlett & & & 10140, 22011 & 96 \\
Sophie and Leo, Chloe, Charlotte & & & 21110, 41020 & 97 \\
& & & 12020, 21111, 41021 & 98 \\
& & & 12021 & 99  \\
\hline
\hline
\textbf{Score averaged over game partners: Single names} & & & \textbf{All $c=3$ sets} & \\
\hline
\hline
Names & Average \% & Set & Set & Average \% \\
\hline
Adam Amy Hannah Harry Jack Max Mia Muhammad Nathan Ruby Ryan & 25.75 & 00000 & 00000 & 25.75  \\
Daisy Isaac & 26.76 & 00001 & 00001 & 26.76 \\
Abigail Dylan Liam Lucy & 35.75 & 10000 & 00002 & 27.76 \\
Isla Lucas & 36.76 & 10001 & 22211 & 34.99  \\
Archie Benjamin Edward Emma Erin Ethan Freya Grace Henry Jake Jamie Jayden Katie Matthew Megan & 37.90 & 10001 & 22210 & 35.57  \\
James Jasmine Maisie Sienna Summer & 38.96 & 00011 & 01222 & 35.74  \\
Jessica & 40.00 & 00012 & 10000 & 35.75  \\
Oliver & 44.55 & 11110 & 12220 & 36.22  \\
Callum Layla Lily William &44.85 &  20000 & 10001 & 36.76 \\
Alfie Alice Amelia Charlie Daniel Emily Finley Lacey Leah Lexy Luke Riley Tyler	& 47.91 & 10010 & 10002 & 37.76  \\
Lewis Samuel Scarlett & 48.97 & 10011 & 22220 & 37.82 \\
\hline
Ellie & 63.53 & 20020 & 02001 & 68.18 \\
Eva Harvey & 63.53 & 00110 & 00201 & 68.39 \\
Imogen Mohammed Theo & 63.82 & 01010 & 11020 & 68.39 \\
George Phoebe & 64.45 & 01020 & 20111 & 68.39 \\
Joseph Sophie & 64.97 & 01011 & 01101 & 68.86 \\
Charlotte Chloe Leo & 65.00 & 11010 & 21011 & 68.86 \\
Olivia & 65.10 & 11100 & 00200 & 69.76 \\
Holly Lola Molly & 65.16 & 21000 & 01102 & 69.76 \\
Evie & 66.69 & 00120 & 20110 & 69.76 \\
Connor & 67.42 & 02000 & 21012 & 69.76 \\
\hline
\hline
\end{tabular}
\caption{Pairs and single names from the UK top 50 boys' and girls' names that score highest and lowest playing `LOVES'. Top section considers pairs of names (complete starting strings); pairs of top-50 names scoring $99\%$, and number strings arising from top-50 pairs that score $\geq 95\%$, are shown. Bottom section considers individual names, and their expected score when playing against all members of the top-50 names lists. The lowest- and highest-scoring ten sets corresponding to top-50 names (left) and from all sets with $c=3$ (right) are shown.}
\label{highscoretable}
\end{table}

\section*{Discussion}
We have explored the behaviour of a popular number game and observed that a simple iterated rule acting on strings of digits produces rich behaviour. Visualising the progress of a game as a trajectory in a space defined by the length of a number string $w$ and the total magnitude of its digits $m$ allows us to identify regions where convergence and divergence is likely, and to characterise the average dynamics underlying movement throughout this space. Based on this picture of the game as a set of steps in $(w,m)$-space, we have found preliminary results for the probability with which the process corresponding to a given string will terminate, and the dynamics through which this termination or divergence comes about. 

Using the most common childrens' names in the UK, we have informally explored the ways in which high and convergent scores can be obtained. For both individual names and name pairs, pairs of `O's are a high-scoring motif; for name pairs, high counts of `L's can also lead to high scores. For individual names, an intermediate number of `LOVES' letters confers higher success than cases where all or none are present. Names with no `LOVES' letters consistently score lowest. The reason for the trend rewarding individual names with intermediate $m$ is unclear, and contrasts with the absence of straightforward predictors of success across complete starting strings (Fig. \ref{figoutcomes}). The specific structure of letter counts among UK names (Figs. \ref{fig2} and \ref{fig2-5}) likely plays a role in determining this optimum.
 
Many open questions remain:

\begin{enumerate}
\item What features of a terminating string predict the magnitude and value of the final result? What gives rise to the complicated structure describing which final results arise more commonly than others?
\item Eqn. \ref{theoryline} is rather empirically derived. More rigorously, what are the parameters describing the transition from highly-likely termination to highly-likely divergence as $m$ increases for a given $w$?
\item What features of a terminating string predict the decay length of the resulting chain?
\item What are the effects of playing the game in a different base?
\item What patterns result when identical letters in the game word enforce correlations in the initial string (for example, `ELSKER', constraining the first and fifth digits to be identical)? 
\end{enumerate}

We hope that this informal study of a popular number game has illustrated its rich behaviour to former and new players alike. We suggest that exploring these results, and extending them, may be an engaging pedagogical route for children and students learning about open-ended mathematical investigation and simulation design \cite{klawe1995classroom}.

\section*{Acknowledgements}
The author wishes to thank C. Brockington, M. Celuzza, and I. Blaikie for suggesting this investigation, and E. R\o yrvik for inspiring it.

\bibliographystyle{unsrt}
\bibliography{loves}

\section*{Appendix}

\subsection*{Number of steps to termination, and chain properties}
\begin{figure}
\centering
\includegraphics[width=0.6\linewidth]{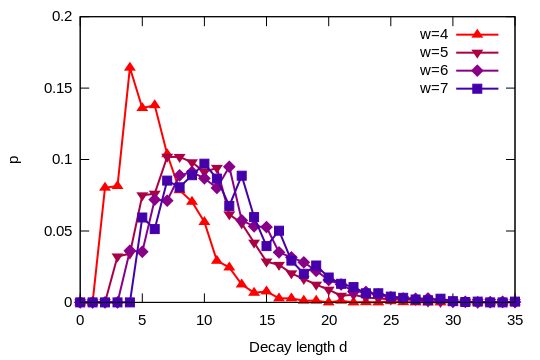}
\caption{Decay lengths (number of iterations required for a string to terminate) of the $c=10$ set for different $w$.}
\label{fig7}
\end{figure}

How long must players keep doing sums before they arrive at their final result? To explore this, we define the `decay length' $d$ as the number of iterations before the game terminates. If a string leads to a looping or divergent outcome, $d = \infty$. Fig. \ref{fig7} shows the distribution of decay lengths for strings in $c=10$ for various $w$ that terminate. The structure of the game places the restriction that $d \geq w-2$: we cannot terminate in fewer steps because we are limited to a maximum string length decrease of 1 at each step. The decay length distribution is moderately skewed, with some strings taking a considerable time to converge (for example, $\{4,0, 3, 9, 0\}$ requires 35 iterations to reach $\{2, 2\}$). A noticeable feature of games displaying this slow decay is the presence of $\{1, 0\}$, $\{1, 1\}$ and $\{1, 2\}$ motifs, often in conjunction with 4s, forming medium-magnitude structures which then combine to exceed 9 and replenish the $\{1, [low]\}$ motifs. 

\begin{figure}
\centering
\includegraphics[width=0.6\linewidth]{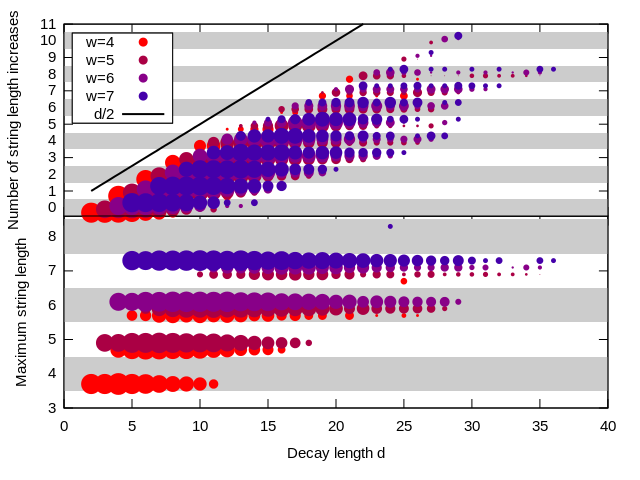}
\caption{(top) The number of string length increases throughout a game, and (bottom) the maximum string length reached during a game, against the decay length (number of iterations required for a string to terminate) of the $c=10$ set for different $w$. Points for different $w$ are offset vertically for clarity.}
\label{figchains}
\end{figure}

What are the dynamics of the chains of strings in these games with long decay lengths? We first explore the maximum string length reached in a chain, as a function of the decay length (Fig. \ref{figchains}). Interestingly, very few terminating chains ever contain strings that are much longer than the initial string -- suggesting that if the initial string length is exceeded by some amount, divergence is guaranteed. Fig. \ref{figchains} shows that terminating chains from $w=4$ only very rarely reach length 7; no terminating chains for $w<7$ reach length 8, and it is very rare for a terminating $w=7$ chain to reach length 8 (and terminating $w=7$ chains are themselves rare). The expansion of a $w=7$ string to a length-8 step occurs in the first iteration of the string $\{1, 0, 0, 0, 1, 9, 1\}$ ($\rightarrow \{1,0, 0, 1, 1, 0, 1, 0\}$, before terminating at $\{2, 5\}$ after 24 iterations), also illustrating the aforementioned appearance of $\{1, [low]\}$ motifs in slow-decaying chains. 

We can also attempt to characterise the dynamics of terminating chains by considering the number of steps in a chain that lead to an increase in string length (Fig. \ref{figchains}). Here, the trend is similar across different $w$ (once the longer decay lengths at higher $w$ are taken into account): the number of increasing steps always falls substantially below half of the total number of steps, suggesting that if the decay chain exceeds a certain ratio of length-increasing steps, divergence is guaranteed.

\subsection*{Critical $m$ for a given $w$}

Fig. \ref{figsteps} allows us to build a heuristic estimate for properties of a trajectory that avoids non-convergence. First note that the vast majority of $w=3$ strings converge; the only exceptions are the $\{9, 9, x\}$ strings, which loop if $1 \leq x \leq 9$. We will neglect these strings and assume that $w = 4$ is the lowest $w$ where a substantial number of strings are non-convergent.

We will consider a trajectory where $w$ decreases from some starting value to $w=2$, without ever increasing. This picture is overly restrictive: we have seen (Fig. \ref{figchains}) that trajectories that experience length increases can still converge. However, for this approximation, we will assume that no length-increasing steps occur ($n=0$). In this case, the maximum possible value of $m$ after an iteration of the game rule is (from Eqn. \ref{mupdate}):

\begin{equation}
m_{new} = m' - 9n = m' = s_1 + s_w + 2\sum_{i=2}^{w-1} s_i. \label{iterate}
\end{equation}

What are the bounds on $m_{new}$? If the current string has $m \leq 18$, it is straightforward to see that $m_{new}$ is minimised when $s_1 + s_w = m$ and $s_i = 0$ for all other $i$ (concentrating the magnitude in the edge digits, which are not doubled by Eqn. \ref{iterate}). This situation gives $m_{new} = m$. 

The maximum $m_{new}$, if $w \geq 4$ and current $m \leq 18$, can be seen to be acheived when the magnitude is concentrated in the sum of Eqn. \ref{iterate}, and thus doubled. This situation is acheived when $s_1 = s_w = 0$ and $\sum_{i=2}^{w-1} s_i = m$, giving $m_{new} = 2m$. If $w = 3$ and current $m \leq 9$, the same reasoning holds. If $w = 3$ and $9 < m \leq 18$, $m_{new}$ is maximised by an arrangement where $\sum_{i=2}^{w-1} = 9$ and $s_1 + s_w = m-9$, giving $m_{new} = m+9$.

Hence, if $m \leq 18$, $m_{new}$ is always bounded from below by $m$ and by above by $2m$ (the $m+9$ bound is enclosed by the $2m$ bound as we have enforced $m > 9$ in this case). Further detailed progress can be made by considering arrangements of digits that induce different values of $n \not= 0$. We will work with the crude approximation that, along a series of non-length-increasing steps, an expected step in $m$ falls in the centre of these bounds. Under this approximation, $m \rightarrow \frac{3}{2} m$.

We have assumed that $w=4$ is the lowest $w$ with a critical value of $m$ above which non-convergence is likely. Considering the repeated action of approximated steps above from some starting $w$, we see that to avoid exceeding a given critical value $m_c$ for $w = 4$, we require $m < m_c / (3/2)^{w-4}$. 

What value of $m_c$ can be chosen to avoid divergence? Consider the maximum $m$ we can achieve for $w = 4$ with $n = 0$ (no length-increasing steps). This is achieved by strings that, upon iteration, give $\{9, 9, 9\}$ -- for example, $\{4, 5, 4, 5\}$ (among others). The maximum $m$ for which $n=0$ at $w=4$ is thus readily seen to be $18$, and we set $m_c = 18$ to follow this line of reasoning. We see in Fig. \ref{figsteps} that $m^* = 18 / (3/2)^{w-4}$ gives a reasonable theoretical prediction of the position of the transitional $m$ for a given $w$.

As something of an aside, we can also place bounds on the minimum magnitude required to support a given value of $n$. It will be observed that the lowest-magnitude strings of adjacent pairs that give rise to $n$ sums $\geq 10$ are $\{1, 9\}$ (among others for $n = 1$); $\{1, 9, 1\}$ ($n = 2$); $\{1, 9, 1, 9\}$ ($n = 3$); $\{1, 9, 1, 9, 1\}$ ($n = 4$); etc. The minimum $m$ required to support a given $n$ is thus:

\begin{equation}
m_{min}(n) = 10 \frac{n+1}{2} + (1 - n \,\mbox{mod}\, 2)
\end{equation}

\subsection*{Backtracking}
Simple algebra shows that, if no steps occur that increase string length, the digits in a decay chain starting from $w=5$ and ending at $w=2$ are given by the expressions in Table \ref{table2}.

A similar set of expressions could be derived in a large number of ways, by choosing the unknown values $a_1, a_2, a_3$ to lie not at the start of each string but at other positions within each. Clearly, to give a valid decay chain with no length-expanding transitions, all the elements in Table \ref{table2} must lie between 0 and 9. The resulting set of inequalities restricting choices of $a_1, a_2,$ and $a_3$ is complicated and, while straightforwardly addressable numerically with integer programming, has prevented further intuitive progress. 

\begin{table}
\scriptsize
\begin{tabular}{|p{3cm}|p{3cm}|p{3cm}|p{3cm}|p{3cm}|}
\hline
$s_1$ & $s_2$ & & & \\
\hline
$a_1$ & $s_1 - a_1$ & $s_2 - s_1 + a_1$ & & \\
\hline
$a_2$ & $a_1 - a_2$ & $s_1 - 2a_1 + a_2$ & $s_2 - 2s_1 + 3 a_1 - a_2$ & \\
\hline
$a_3$ & $a_2 - a_3$ & $a_1 - 2 a_2 + a_3$ & $s_1 - 3a_1 + 3a_2 - a_3$ & $s_2 - 3 s_1 + 6 a_1 - 4 a_2 + a_3$ \\
\hline
\end{tabular}
\caption{Backtracked values of a $w=5$ string giving a final result $\{s_1, s_2\}$ with the assumption of no length-increasing steps. Table rows run backwards through iterations of the game rules; at each backwards step, a new variable $a_i$ is introduced. All elements in this table must lie between 0 and 9 inclusive for this decay chain to be valid.}
\label{table2} 
\end{table}

\end{document}